\documentclass[11 pt, psamsfonts]{amsart}
\usepackage{amsmath}
\usepackage{amsthm}
\usepackage{amsfonts}
\usepackage{amssymb}
\usepackage{amscd,epsf, latexsym}
\usepackage[dvips]{graphicx}

\usepackage{eucal}

\newcommand{\B}{\mathcal{B}}

\newcommand{\PP}{\mathcal{P}}

\newcommand{\lan}{\langle}
\newcommand{\ra}{\rangle}

\newcommand{\tr}{{\rm tr}}
\newcommand{\Ker}{{\rm Ker}}

\newcommand{\Ind}{{\rm Ind}}
\newcommand{\Irr}{{\rm Irr}}
\newcommand{\Res}{{\rm Res}}

\newcommand{\Aut}{{\rm Aut}}
\newcommand{\Inn}{{\rm Inn}}

\newcommand{\Z}{\mathbb{Z}}

\newcommand{\C}{\mathbb{C}}

\newcommand{\bE}{\mathbf{E}}
\newcommand{\la}{{\lambda}}

\newcommand{\al}{\alpha}

\theoremstyle{plain}
\newtheorem{theorem}{Theorem}[section]
\newtheorem{lemma}[theorem]{Lemma}

\newtheorem{prop}[theorem]{Proposition}
\newtheorem{cor}[theorem]{Corollary}
\theoremstyle{definition}
\newtheorem{definition}[theorem]{Definition}
\theoremstyle{remark}

\newtheorem{rmk}[theorem]{Remark}

\title{Extraspecial 2-groups and images of braid group representations}

\author{Jennifer Franko}
\email{jefranko@indiana.edu}
\address{Department of Mathematics\\Indiana University \\Bloomington, IN 47405\\U.S.A.}

\author{Eric C. Rowell}
\email{errowell@indiana.edu}
\address{Department of Mathematics\\Indiana University \\Bloomington, IN 47405\\U.S.A.}

\author{Zhenghan Wang}
\email{zhewang@indiana.edu}
\address{Department of Mathematics\\Indiana University \\Bloomington, IN 47405\\U.S.A.}
\thanks{The first author is supported in part by NSF grant EIA 0130388.
\\  The second author is supported in part by an NSF VIGRE grant.
\\  The third author is partially supported by NSF grants EIA 0130388 and FRG 0354772.}

\begin{document}

\begin{abstract} We investigate a family of (reducible) representations of the braid groups $\B_n$ corresponding to a specific solution to the Yang-Baxter equation.  The images of
$\B_n$ under these representations are finite groups, and we identify them precisely as extensions of extra-special 2-groups.  The decompositions of the representations into their irreducible constituents are determined, which allows us to relate them to the well-known Jones representations of $\B_n$ factoring over Temperley-Lieb algebras and the corresponding link invariants.
\end{abstract}

\maketitle

\section{Introduction}

Representations of Artin's braid groups $\B_n$ are of great importance
to mathematicians \cite{Bir}, and physicists recently \cite{W}.  Certain representations of the braid
groups have been proposed as the fractional
statistics of anyons \cite{W}, and used in the topological models for
quantum computing \cite{FKLW}.  Therefore it is interesting to identify the
images of such braid group representations.  In this paper
we analyze a particular representation of the braid groups
afforded by a unitary solution
of the braid relation, i.e. a flipped R-matrix
$R=\frac{1}{\sqrt{2}}\begin{pmatrix}1 & 0 & 0 & 1\\0 & 1& -1 & 0\\0 & 1&1& 0\\-1 & 0 & 0
&1\end{pmatrix}$
that satisfies the \emph{Yang-Baxter equation}
\begin{enumerate}
\item[(YBE)] $\quad(R\otimes I_2)(I_2\otimes R)(R\otimes I_2)=(I_2\otimes R)(R\otimes I_2)(I_2\otimes R)$
\end{enumerate} where $I_2$ is the $2\times 2$ identity matrix.
All solutions to the YBE of the form $R: V\otimes V\rightarrow V\otimes V$
with $V$ 2-dimensional have been listed
in \cite{H}.  Dye \cite{D} found all unitary solutions of this form to the braid relations
based on this list.  The  importance of this particular braid
operator $R$ was pointed out in the work of Kauffman and Lomonaco
\cite{KauffLo}, and the connection of $R$ with quantum computing was explored
there which is another reason for our interest.

As is well-known,
any (invertible) matrix satisfying the YBE gives rise to representations of $\B_n$ for any $n$.
The  representation $(\pi_n,(\C^2)^{\otimes n})$ corresponding to the matrix $R$ above
is unitary and defined as follows:
$$\pi_n(\sigma_i)=I_2^{\otimes i-1}\otimes R\otimes I_2^{\otimes n-i-1},$$
where $\sigma_i$ is the $i$-th braid generator.
The images of the braid groups $\B_n$ under this representation
are finite groups,
and the image matrix of each braid generator $\sigma_i$ has only two distinct eigenvalues.
It follows that the image group of an irreducible constituent of $\pi_n$ is
generated by the conjugacy class of a braid generator with two distinct
eigenvalues whose ratio is not -1, \emph{i.e.} has the so-called
2-eigenvalue property defined in \cite{FLW}.  Such representations
are completely classified \cite{FLW}, so in principle the image groups of the
irreducible constituents of $\pi_n$ can be
identified by using the complete list in \cite{FLW} Theorem 1.6.  But as
we will see that $\pi_n$ is reducible, hence first we need to
find the irreducible constituents of $\pi_n$;
then we need to distinguish a few different cases in the complete
list for the images of the irreducible constituents, so instead
we choose to solve the problem in an elementary and self-contained way.
We decompose these representations $\pi_n$ (for all $n$)
into their irreducible constituents and describe the images
of $\B_n$ under $\pi_n$ as abstract groups.
We find that the images of the {\it pure} braid
groups are (nearly) extra-special 2-groups $\bE_{n-1}^{-1}$.  The images of
the full braid groups $\B_n$ are extensions of the (nearly) extra-special 2-groups
$\bE_{n-1}^{-1}$ by the symmetric groups $S_n$, and the restrictions of the representations $\pi_n$ to the subgroup of pure braids are isotypic
copies of the odd representations of $\bE_{n-1}^{-1}$.

As already discussed in \cite{KauffLo} we can define link invariants
using the representations $\pi_n$.
By observing that $\pi_n$ is related to the Jones
representation of the braid groups at the 4-th root of unity, we
improve slightly some earlier results of Jones about the images of the
Jones representation of the braid groups at the 4-th root of unity \cite{Jones86}.
As a consequence we point out that the resulting link invariants are essentially
the Jones polynomial at a 4-th root of unity, hence really the
Arf invariant of a link (see references in \cite{Jones}).  The slight improvement of Jones's
result comes from two subtle points about the Jones representations.
Firstly, in the Jones representation of the braid group, there is some
freedom in choosing phases so it is
convenient to
 state the results projectively, \emph{i.e.} modulo scalars, while not
 losing any significance mathematically.  We choose to work out
the images in full generality (as opposed to projectively) as this is desirable in physics
for the applications to the fractional statistics of quantum Hall
fluid \cite{Read}.  This changes the images of the pure braid groups
from the elementary abelian groups $\Z_2^{n-1}$ to the
(nearly) extra-special 2-groups $\bE_{n-1}^1$.
Secondly, when the number of strands of the braid
groups is even, there are two irreducible sectors of the Jones
representation \cite{Jones86}.  Jones found the projective images for each
sector, but we determine the images of the two sectors together.
This brings up a subtlety about the centers of the (nearly) extraspecial
2-groups in those cases, which disappears when the two irreducible
sectors are treated separately, and projectively.

These results lead to several questions for future research currently being worked out by the authors.
What are the closed images of the braid groups under the representations afforded by
the other $R$-matrices listed in \cite{D} and
what are the associated link invariants?  What are the other extraspecial $p$-groups
that appear as homomorphic images of the pure braid groups?

\section{Preliminaries}
\subsection{Definitions and computations}\label{defs}

\begin{definition}
Artin's braid group $\B_n$ on $n$ strands has presentation in generators $\sigma_1,\ldots,\sigma_{n-1}$
satisfying relations:
\begin{enumerate}
\item[(B1)] $\sigma_i\sigma_j=\sigma_j\sigma_i$ if $|i-j|\geq 2$.
\item[(B2)] $\sigma_i\sigma_{i+1}\sigma_i=\sigma_{i+1}\sigma_i\sigma_{i+1}$ for $1\leq i\leq n-2$
\end{enumerate}
\end{definition}

For tensor products of matrices we use the convention ``left into right," that is,
if $X=\begin{pmatrix}w & x \\ y & z\end{pmatrix}$ and $A=\begin{pmatrix}a & b \\ c & d\end{pmatrix}$
then $X\otimes A=\begin{pmatrix} aX & bX \\ cX & dX\end{pmatrix}$.
Various matrices and quantities will be needed throughout, so we define them here:
\begin{enumerate}
\item $I_m$ is the $m\times m$ identity matrix.
\item $R=\frac{1}{\sqrt{2}}\begin{pmatrix}1 & 0 & 0 & 1\\0 & 1& -1 & 0\\0 & 1&1& 0\\-1 & 0 & 0 &1\end{pmatrix}$
\item $s=\begin{pmatrix}0 & 1\\ -1 &0\end{pmatrix}$
\item $\sigma_x=\begin{pmatrix}0 & 1\\ 1 &0\end{pmatrix}$
\item $P_s=\begin{pmatrix}1 & \sqrt{-1}\\ \sqrt{-1} &1\end{pmatrix}$
\item $P_{\sigma_x}=\begin{pmatrix}1 & -1\\
1 & 1\end{pmatrix}$
\item $P_n=(P_s\otimes P_{\sigma_x})^{\otimes \lfloor n/2\rfloor}\otimes I_2^{\otimes (n-2\lfloor n/2\rfloor)}$ where $\lfloor a\rfloor$ is the integer part of $a$.
\item $\sigma_z=\begin{pmatrix}1 & 0\\ 0 &-1\end{pmatrix}$
 \item $g_i=\pi_n((\sigma_i)^2)=I_2^{\otimes (i-1)}\otimes R^2\otimes I_2^{\otimes (n-i-1)}$
 (observe we ignore the dependence of $g_i$ on $n$; the value of $n$ will always be clear from the context).
\item $\zeta=\frac{1}{\sqrt{2}}(1+\sqrt{-1})$.
 \item $d=\begin{pmatrix}\zeta & 0\\ 0 & \overline{\zeta}\end{pmatrix}$
 \item $D=\begin{pmatrix}\zeta & 0&0&0\\ 0 & \overline{\zeta}&0&0\\
 0&0&\overline{\zeta}&0\\
 0&&0&\zeta\end{pmatrix}$
 \item $M=\frac{1}{\sqrt{2}}\begin{pmatrix}1 & 1\\ -1 & 1\end{pmatrix}$
 \end{enumerate}
We will also need a few simple computations involving these matrices so we record them in the following:
\begin{lemma}\label{comps}The matrices defined above satisfy:
\begin{enumerate}

\item[(a)]
$R^2=s\otimes\sigma_x$
\item[(b)] $R=\frac{1}{\sqrt{2}}(R^2+I_4)$, $R^{-1}=\frac{1}{\sqrt{2}}(R^{-2}+I_4)$
\item[(c)] $(R^2\otimes I_2)(I_2\otimes R^2)=-(I_2\otimes R^2)(R^2\otimes I_2)$
\item[(d)] $g_ig_{i+1}=-g_{i+1}g_i$
\item[(e)] $(R^{-1}\otimes I_2)(I_2\otimes R^2)(R\otimes I_2)=(I_2\otimes R^2)(R^2\otimes I_2)$
\item[(f)] $\pi_n(\sigma_i^{-1})g_{i\pm 1}\pi_n(\sigma_i)=g_{i\pm 1}g_i$
\item[(g)] $g_ig_j=g_jg_i$ and $\pi_n(\sigma_i)g_j=g_j\pi_n(\sigma_i)$ if $|i-j|\geq 2$
\item[(h)] $R^4=-I_4$, $(g_i)^2=-I_{2^n}$.
\item[(i)] $(P_s)^{-1}sP_s=\sqrt{-1}\sigma_z$, and $(P_s)^{-1}\sigma_xP_s=\sigma_x$
\item[(j)] $(P_{\sigma_x})^{-1}\sigma_xP_{\sigma_x}=\sigma_z$, $(P_{\sigma_x})^{-1}sP_{\sigma_x}=s$
\item[(k)] $(P_n)^{-1} g_{2i+1}P_n=\sqrt{-1}(I_2^{\otimes 2i}\otimes \sigma_z^{\otimes 2}\otimes I_2^{\otimes n-2i-2})$, and $(P_n)^{-1} g_{2i}P_n=g_{2i}$.
\end{enumerate}\end{lemma}
\begin{proof}
The first assertions (a) and (b) are straightforward computations.
Having checked that $s$ and $\sigma_x$ anti-commute (c) follows, and
(d) is immediate from (c). Using (b) and the observation $R^{-2}=-R^2$
we express the left-hand side of the equality in (e) in terms of $R^2$ and
then use (c) to derive the right-hand side.  Assertion (f) is immediate from (e).
Assertion (g) is a consequence of the ``far commutation" relations satisfied by the braid group,
and (h) follows from (b) and the definition of $g_i$.
The matrix $P_s$ (resp. $P_{\sigma_x}$) is a change of bases matrix that diagonalizes $s$
(resp. $\sigma_x$) and commutes with $\sigma_x$ (resp. $s$).
This is the statement (j), and (k) follows directly from this fact and the definition of $P_n$.
\end{proof}

\subsection{Restriction to $\PP_n$}
The homomorphism $\B_n$ to the symmetric group on $n$ letters $S_n$ given by $\sigma_i\rightarrow (i,i+1)$
has kernel $\PP_n$ the so-called pure braid group.  $\PP_n$ is generated by all conjugates of the
squares of the generators of $\B_n$: $(\sigma_i)^2$.  Actually a more economical presentation
of $\PP_n$ can be found (see e.g. \cite{Bir}), but we shall not need it here.
To exploit this relationship between $\B_n$ and $\PP_n$ we shall restrict $\pi_n$ to the subgroup $\PP_n$.
For convenience of notation we introduce the following notation:
\begin{definition}
$H_n:=\pi_n(\PP_n)$ and $G_n:=\pi_n(\B_n)$
\end{definition}
 We can describe $H_n$ very succinctly:
\begin{lemma}
$H_n$ is generated by $g_1,\ldots,g_{n-1}$.
\end{lemma}
\begin{proof}  Observe that $H_n$ is generated by all conjugates of $g_i$, so that $H_n$ is the smallest normal subgroup of $G_n$ containing the subgroup $\lan g_1,\ldots,g_{n-1}\ra$ generated by the $g_i$.  But by Lemma \ref{comps}(f),(g) $\lan g_1,\ldots,g_{n-1}\ra$ is normal in $G_n$ so $\lan g_1,\ldots,g_{n-1}\ra=H_n$.
\end{proof}
\begin{rmk}
Combining this with Lemma \ref{comps}(a) we have a very powerful tool for studying the representation $\pi_n$ of $\B_n$.  After decomposing the representation $\pi_n$ restricted to $\PP_n$ into its irreducible components and computing the corresponding images of the $(\sigma_i)^2$, we can immediately determine the decomposition of the images of the $\sigma_i$ under $\pi_n$ as $\pi_n(\sigma_i)=\frac{1}{\sqrt{2}}(g_i+I_{2^n})$.
\end{rmk}

Once we understand $H_n$ as an abstract group and decompose its defining representation
(as it is presented to us as a matrix group) we will need to consider the group $G_n/H_n$.
We can immediately see that  $G_n/H_n$ is a homomorphic image of $S_n$ as $\pi_n$
induces a surjective homomorphism
$\overline{\pi}_n: \B_n/\PP_n\rightarrow G_n/H_n$ and $\B_n/\PP_n\cong S_n$.
We would like to know if $\overline{\pi}_n$ is an isomorphism, so we must determine if $\Ker(\overline{\pi}_n)$ is
trivial.  Observing that the kernel of $\overline{\pi}_n$ is (isomorphic to) a normal subgroup of $S_n$
we need only check that the kernel is not $S_n$, $A_n$ or the normal subgroup of $S_4$
isomorphic to $\Z_2\times\Z_2$.  For $n\geq 4$ it is sufficient to check that the element $(12)(34)$
is not in the kernel, while for $n=3$ we should check that $(123)$ is not in the kernel.
Since $H_n$ is a normal subgroup of $G_n$, we also have a homomorphism $\vartheta: G_n\rightarrow \Aut(H_n)$
where $\Aut(H_n)$ is the automorphism group of $H_n$ and $\vartheta(s)$ is conjugation by $s\in G_n$.
Restricting to $H_n$ we see that $\vartheta(H_n)=\Inn(H_n)\subset\Aut(H_n)$ the normal subgroup of
inner automorphisms of $H_n$ and so we have the induced homomorphism
$\overline{\vartheta}:G_n/H_n\rightarrow \Aut(H_n)/\Inn(H_n)$.
Since $\Ker(\overline{\pi}_n)\subset\Ker(\overline{\vartheta}\circ\overline{\pi}_n)$
if we can show the composition has trivial kernel then $\overline{\pi}_n$ must be an isomorphism.
By Lemma \ref{comps}(d) the generators $g_i$ of $H_n$ commute or anti-commute,
so the elements of $\Inn(H_n)$ act by sign changes.
So if we can show that the automorphisms corresponding to $(12)(34)$ (for $n\geq 4$) and $(123)$ are not
simply sign changes, we will have shown that $\overline{\pi}_n$ is an isomorphism.
The corresponding elements of $\B_n$ are $(\sigma_1\sigma_3)$ and $(\sigma_2\sigma_1)$ and
we use Lemma \ref{comps}(f) to compute that under
$\overline{\vartheta}\overline{\pi}_n$ the element $(\sigma_1\sigma_3)$ maps $g_2$ to $g_2g_1g_3$,
and $(\sigma_2\sigma_1)$ maps $g_2$ to $g_2g_1g_2=g_1$. We check directly that $g_1g_3\not=\pm 1$ using Lemma \ref{comps}(a),(i) and (j), so $\overline{\pi}_n$ is an isomorphism for $n\geq 3$.  In the case $n=2$ we see that $G_2$ is the group generated by the matrix $R$ which is isomorphic to $\Z_4$, so combining, we have:
\begin{theorem}\label{quotient}
We have an exact sequence:
$$1\rightarrow H_n\stackrel{\subset}{\rightarrow}G_n\stackrel{\vartheta}{\rightarrow} S_n\rightarrow 1$$
for all $n\geq 2$.  In other words, $G_n$ is an extension of $H_n$ by $S_n$.
\end{theorem}

\section{Extraspecial 2-groups and related groups}

\begin{definition}
The group $\bE_m^\nu$ is the abstract group generated by $$x_1,\ldots,x_{m}$$
with relations:
\begin{eqnarray}
&&x_i^2=\nu, \quad 1\leq i\leq m\\
&&x_ix_j=x_jx_i, \quad |i-j|\geq 2\\
&&x_{i+1}x_i=-x_ix_{i+1}, \quad 1\leq i\leq m
\end{eqnarray}
where $-1$ is an order two central element, and $\nu=\pm 1$.
\end{definition}
These groups appear classically and have important connections with Clifford algebras.  The case $\nu=-1$
appears in Exercise 3.9 in the text by Fulton and Harris \cite{FH}, and
other cases appeared in \cite{Griess}.  The necessary facts about these groups are found in various places,
but are elementary so we reprove them here for the reader's convenience.

\subsection{Properties of $\bE_m^\nu$}
Any element in $\bE_m^\nu$ can be expressed in the normal form: $\pm x_1^{\al_1}\cdots x_{m}^{\al_{m}}$
where $\al_i\in\Z_2$.  The following lemma will show that it is unique.
\begin{lemma}\label{centerES}
Denote by $Z(\bE_m^\nu)$ the center of $\bE_m^\nu$.  We have:
\begin{enumerate}
\item[(a)]
$Z(\bE_m^\nu)=\begin{cases}\{\pm 1\} &  \textrm{$m$ even}\\
\{\pm 1, \pm x_1x_3\cdots x_m\} & \textrm{$m$ odd}
\end{cases}$
\item[(b)] $\bE_m^\nu/\{\pm 1\}\cong (\Z_2)^{m}$
\item[(c)]  Any $x\in\bE_m^\nu\setminus Z(\bE_m^\nu)$ is conjugate to $-x$.
\item[(d)] Any nontrivial normal subgroup of $\bE_m^\nu$ intersects $Z(\bE_m^\nu)$ nontrivially.
\item[(e)] For $m=2k-1$ odd, $Z(\bE_{2k-1}^\nu)\cong\begin{cases}\Z_2\times\Z_2 & \textrm{if $\nu=1$ or $k$ even}\\
\Z_4 & \textrm{if $\nu=-1$ and $k$ odd}
\end{cases}$
\item[(f)] The normal form $\pm x_1^{\al_1}\cdots x_m^{\al_m}$ is unique.
\end{enumerate}
\end{lemma}
\begin{proof}
Using the above-mentioned normal form we may assume, without loss of generality,
that $z=x_1^{\al_1}\cdots x_{2k}^{\al_{m}}\in Z(\bE_m^\nu)$ since if $z$ is central, so is $-z$.
By the commutation/anti-commutation relations in $\bE_m^\nu$
we have $x_iz=(-1)^{\al_{i-1}+\al_{i+1}}zx_i=zx_i$ for all $i$ where we take $\al_{-1}=\al_{m+1}=0$.
Thus we get the system of equations over $\Z_2$:
\begin{eqnarray*}
&&\al_2=0\\
&&\al_{m-1}=0\\
&&\al_i+\al_{i+2}=0\pmod{2}, \quad 1\leq i\leq m-2
\end{eqnarray*}
If $m$ is even then the system has only the trivial solution $\boldsymbol{\al}=\mathbf{0}$,
but if $m$ is odd there are two solutions $\mathbf{0}$ and $(1,0,1,\ldots,0,1)$,
that is, all the $\al_{2i}=0$ and $\al_{2i+1}=1$.   Thus we have (a).
It is clear from the relations in $\bE_m^\nu$ that $\bE_m^\nu/\{\pm 1\}$ is presented
by $m$ commuting generators of order 2, \emph{i.e} $(\Z_2)^{m}$.
To prove (c) observe that any non-central element $x\in\bE_m^\nu$ must
anti-commute with some $x_i$.  So (d) follows from (c) as any nontrivial normal subgroup $N$
must either be central or contain $\{x,-x\}$ for some non-central element $x$ so that $-1\in N$ as well.
For (e) we compute the order of the central element $x_1x_3\cdots x_{2k-1}$ and find that it is 2 or 4,
which gives us the two cases.  Assertion (f) follows from a simple counting argument
as $|\bE_m^\nu|=2^{m+1}$ by (b).
\end{proof}

 \begin{definition} A group $G$ of order $2^{m+1}$ is an \textbf{extraspecial 2-group}
 if (see \cite{Griess}):
 \begin{enumerate}
 \item The center $Z(G)$ and the commutator subgroup $G^\prime$ coincide and are isomorphic to $\Z_2$.
 \item $G/Z(G)\cong (\Z_2)^m$.
 \end{enumerate}
 \end{definition}
It is immediate from the anti-commutation relations that the commutator subgroup
of $\bE_m^\nu$ is $\{\pm 1\}$, and for $m=2k$ the other conditions were verified
in Lemma \ref{centerES} so we have:
\begin{prop}
$\bE_{2k}^\nu$ is an extraspecial 2-group.
\end{prop}
\begin{rmk}
Since the group $\bE_{2k+1}^{\nu}$ contain $\bE_{2k}^\nu$, we call the groups $\bE_m^\nu$
\emph{nearly} extraspecial 2-groups for any $m$ (so they include extraspecial 2-groups).  This should not be confused with \emph{almost} extraspecial 2-groups found in the literature which are central products of extraspecial 2-groups with $\Z_4$.  The cases where the center of $\bE_m^\nu$ is isomorphic to $\Z_4$ are almost extraspecial, but when the center is $\Z_2\times\Z_2$ they are not (see \cite{Griess}).
\end{rmk}
\subsection{Representations of $\bE_{m}^{\nu}$}\label{esk}

We wish to construct the irreducible representations of $\bE_{m}^{\nu}$.  There are 4 cases
corresponding to the parity of $m$ and the choice of $\nu$.
For the reader's convenience we recall the following standard facts from the character theory of
finite groups (see any standard text, e.g. \cite{FH}):
\begin{prop}\label{chars}
Let $G$ be a finite group, and $\Irr(G)=\{\chi_i\}_{i\in\mathcal{I}}$ the set of irreducible
characters of $G$, corresponding to irreducible representations $V_i$.
\begin{enumerate}
\item[(a)] $|\Irr(G)|$ is equal to the number of conjugacy classes of elements of $G$.
\item[(b)] $|G|=\sum_\mathcal{I}(\dim V_i)^2$
\item[(c)] For $\chi_i$, $\chi_j\in\Irr(G)$ $\sum_{g\in G}\chi_i(g)\overline{\chi_j(g)}=\begin{cases} 0 & \textrm{if $V_i\not\cong V_j$}\\ |G| & \textrm{if $V_i\cong V_j$}\end{cases}$.
\item[(d)] If $g$ and $h$ are not conjugate then $\sum_\mathcal{I}\chi_i(h)\overline{\chi_i(g)}=0$.
\end{enumerate}
\end{prop}

\subsubsection{$\bE_{2k}^{-1}$}
To determine the number of irreducible representations we count conjugacy classes.  The center $\{\pm 1\}$
gives us two singleton classes, and Lemma \ref{centerES}(c) shows that the non-singleton conjugacy classes
are given by $[\pm x_1^{\al_1}\cdots x_{2k}^{\al_{2k}}]$ for
any $\boldsymbol{\al}\in (\Z_2)^{2k}\setminus\{\mathbf{0}\}$.
So we have $2+(2^{2k}-1)=2^{2k}+1$ inequivalent irreducible representations.
Let $\Irr(\bE_{2k}^{-1})=\{V_1,\cdots,V_{2^{2k}+1}\}$ denote a set of
inequivalent irreducible representations of $\bE_{2k}^{-1}$.
By Lemma \ref{centerES} we can induce 1-dimensional representations of $\bE_{2k}^{-1}$
from any representation of $(\Z)^{2k}$ by letting the center act trivially.
Thus we have $2^{2k}$ 1-dimensional representations (say, $V_2,\ldots, V_{2^{2k}+1}$)
leaving only one representation, $V_1$ to determine.  Using the class equation
$$2^{2k+1}=|\bE_{2k}^{-1}|=(\dim V_1)^2+\sum_2^{2^{2k}+1}(\dim V_i)^2=(\dim V_1)^2+2^{2k}$$
we find that
$\dim V_1=2^k$.  The 1-dimensional representations are equal to their characters so
for $2\leq i\leq 2^{2k}+1$ we have $\chi_i(1)=\chi_i(-1)=1$, and $\chi_i([\pm x_j])=\pm 1$
for all possible choices of sign.
From Proposition \ref{chars}(c),(d) we find that $\chi_1(1)=-\chi_1(-1)=2^k$,
and $\chi_1([\pm x_1^{\al_1}\cdots x_{2k}^{\al_{2k}}])=0$.
We can construct the representation $(\rho_1,V_1)$
as follows (recall the definition of the matrices $s$ and $\sigma_z$ from section \ref{defs}):
\begin{eqnarray*}
&&\rho_1(x_{1})=\sqrt{-1}(\sigma_z\otimes I_2^{\otimes k-1})\\
&&\rho_1(x_{2})=s\otimes I_2^{\otimes k-1}\\
&&\quad \vdots\\
&&\rho_1(x_{2i})=I_2^{\otimes i-1}\otimes s \otimes I_2^{\otimes k-i}\\
&&\rho_1(x_{2i+1})=\sqrt{-1}(I_2^{\otimes i-1}\otimes \sigma_z\otimes \sigma_z\otimes I_2^{\otimes k-i-1})\\
&&\quad \vdots\\
&&\rho_1(x_{2k})=I_2^{\otimes k-1}\otimes s
\end{eqnarray*}
 As $(\sigma_z)^2=I_2$, $s^2=-I_2$ and $\sigma_z s\sigma_z=-s$ we see that this is indeed
 a representation of $\bE_{2k}^{-1}$, and since $tr(s)=tr(\sigma_z)=0$ it follows from the orthogonality
 of characters that this is the irreducible $2^k$-dimensional representation of $\bE_{2k}^{-1}$.

\subsubsection{$\bE_{2k-1}^{-1}$}
We now construct the irreducible representations of $\bE_{2k-1}^{-1}$.
Denote by $z$ the central element $x_1x_3\cdots x_{2k-1}$ for convenience.
Using Lemma \ref{centerES} we find that there are $2^{2k-1}+2$ distinct conjugacy
classes in $\bE_{2k-1}^{-1}$ and therefore we may label the inequivalent classes of
irreducible representations by $\Irr(\bE_{2k-1}^{-1})=\{W_1,\ldots,W_{(2^{2k-1}+2)}\}$.
We get $2^{2k-1}$ distinct 1-dimensional representations from $(\Z_2)^{2k-1}$ by composing
with the projection onto $\bE_{2k}/\{\pm 1\}$, and them by $W_3,\ldots,W_{(2^{2k-1}+2)}$.
We compute their characters $\psi_i$ for $3\leq i\leq 2^{2k-1}+2$ as in the $\bE_{2k}^{-1}$ case
$\psi_i(1)=\psi_i(-1)=1$ and $\psi_i(\pm x_j)=\pm 1$ which determines
their values on all classes (observe that $\psi_i(z)=\psi_i(-z)$ for nontrivial central
elements $\pm z$ for these 1-dimensional representations).
From Proposition \ref{chars}(b) we get $\dim W_1+\dim W_2=2^k$ for the remaining two
irreducible representations.
Since $\dim W_i | 2^{2k}$ we see that in fact, $\dim W_1=\dim W_2=2^{k-1}$.
Using Proposition \ref{chars}(c) we find that the characters $\psi_1$ and $\psi_2$ vanish
on all equivalence classes except for the central classes: $[1]$, $[-1]$, $[z]$ and $[-z]$.
Observing that the restrictions of $W_1$ and $W_2$ to the
subgroup $\bE_{2k-2}^{-1}\subset\bE_{2k-1}^{-1}$ must both be
the unique non-trivial irreducible $2^{k-1}$-dimensional
representation we find that $\psi_1(-1)=\psi_2(-1)=-2^{k-1}$.
Proposition \ref{chars}(c),(d) then implies first that $\psi_1(z)=\psi_2(-z)=-\psi_1(-z)=-\psi_2(z)$,
and then using this and the orthogonality of $\psi_1$ and $\psi_2$ to see that $\|\psi_1(z)\|=2^{k-1}$.
Restricting to $Z(\bE_{2k-1}^{-1})$ and recalling that $\Z_2\times\Z_2$ has only real characters
while the non-trivial characters of $\Z_4$ have pure complex values on its generators
we determine the value of $\psi_1(z)$ up to a choice of sign coming from switching $W_1$ and $W_2$.
For the purpose of simplifying notation later we include
a sign depending on the value of $k\pmod{4}$ and define:
\begin{eqnarray}\label{psi1}
\psi_1(x)=\begin{cases}\pm 2^{k-1} & \textrm{for $x=\pm 1$}\\
\pm(-1)^{(k/2}(2^{k-1}) & \textrm{for $x=\pm z$}\\
0 & \textrm{otherwise}
\end{cases}
\end{eqnarray}
and
\begin{eqnarray}\label{psi2}
\psi_2(x)=\begin{cases}\pm 2^{k-1} & \textrm{for $x=\pm 1$}\\
\mp(-1)^{(k/2}(2^{k-1}) & \textrm{for $x=\pm z$}\\
0 & \textrm{otherwise}
\end{cases}\end{eqnarray}

Next we give explicit matrix realizations of $W_1$ and $W_2$.
Since $Z(\bE_{2k-1}^{-1})$ must act non-trivially (although not necessarily faithfully)
on $W_1$ and $W_2$ we use the inclusion $\bE_{2k-1}^{-1}\subset\bE_{2k}^{-1}$ to observe:
$$\Ind_{\bE_{2k-1}^{-1}}^{\bE_{2k}^{-1}}(W_1)=\Ind_{\bE_{2k-1}^{-1}}^{\bE_{2k}^{-1}}(W_2)=V_1$$
where $V_1$ is the $2^k$-dimensional irreducible representation of $\bE_{2k}^{-1}$ given in \ref{esk}.
Thus by Frobenius reciprocity (and a dimension count) we have that

\begin{eqnarray}\label{resV}
\Res_{\bE_{2k-1}^{-1}}^{\bE_{2k}^{-1}}(V_1)=W_1\oplus W_2
\end{eqnarray}

From this we get explicit realizations $(\la_1,W_1)$ and $(\la_2,W_2)$.
(\emph{N.b.} the only difference of $\la_1$ and $\la_2$ on the generators is
that the image of $x_{2k-1}$ differs in sign.)
\begin{eqnarray*}
&&\la_1(x_{1})=\la_2(x_1)=\sqrt{-1}\sigma_z\otimes I_2^{\otimes k-2}\\
&&\la_1(x_{2})=\la_2(x_2)=s\otimes I_2^{\otimes k-2}\\
&&\quad \vdots\\
&&\la_1(x_{2i})=\la_2(x_{2i})=I_2^{\otimes i-1}\otimes s \otimes I_2^{\otimes k-i-1}\\
&&\la_1(x_{2i+1})=\la_2(x_{2i+1})=\sqrt{-1}I_2^{\otimes i-1}\otimes \sigma_z\otimes \sigma_z\otimes I_2^{\otimes k-i-2}\\
&&\quad \vdots\\
&&\la_1(x_{2k-2})=\la_2(x_{2k-2})=I_2^{\otimes k-2}\otimes s\\
&&\la_1(x_{2k-1})=-\la_2(x_{2k-1})=\sqrt{-1}I_2^{\otimes k-2}\otimes \sigma_z
\end{eqnarray*}
One easily checks that these indeed define irreducible representations of $\bE_{2k-1}^{-1}$
just as in the $m=2k$ case.  It is perhaps worth computing the traces of the images
of the central element $z$ under $\la_1$ and $\la_2$.  We have:
$\la_1(z)=-\la_2(z)=(\sqrt{-1})^k((\sigma_z)^2\otimes\cdots\otimes (\sigma_z)^2)=(\sqrt{-1})^kI_{2^{k-1}}$
so that:
$$\tr(\la_1(z))=-\tr(\la_2(z))=\begin{cases}2^{k-1} & \textrm{if $k\equiv 0\pmod{4}$}\\
-2^{k-1} & \textrm{if $k\equiv 2\pmod{4}$}\\
\sqrt{-1}(2^{k-1}) & \textrm{if $k\equiv 1\pmod{4}$}\\-\sqrt{-1}(2^{k-1}) & \textrm{if $k\equiv 3\pmod{4}$}
\end{cases}$$
The traces of the images of $\pm 1$ are also easily computed, and
comparing these values with the above formulas \ref{psi1} and \ref{psi2}, we check that the
characters of $\la_1$ and $\la_2$ are $\psi_1$ and $\psi_2$ respectively.

\subsubsection{$\bE_{m}^{1}$}\label{nuis1reps}
Suppose that $(\rho,V)$ is any representation of $\bE_{m}^{-1}$ defined on generators $\rho(x_i)=A_i$
for some set of matrices $\{A_i\}_{1\leq i\leq m}$.  Denote by $x_1^\prime,\ldots,x_{m}^\prime$
the generators of $\bE_{m}^{1}$ and define $\rho^\prime(x_i^\prime)=\sqrt{-1}A_i$.
Then since $(A_i)^2=-Id_V$ we have $(\rho^\prime(x_i^\prime))^2=Id_V$ and $(\rho^\prime,V)$
defines a representation of $\bE_m^{1}$ (observe that the commutation/anti-commutation
relations are homogeneous and hence also satisfied).
Obviously this process is reversible, so that all representations of $\bE_m^{1}$ are obtained in this way.
If we define representations $\la_1^\prime$ and $\la_2^\prime$ of $\bE_{2k-1}^{1}$
corresponding to the two $2^{k-1}$-dimensional representations of $\bE_{2k-1}^{-1}$
then we find that the characters $\psi^\prime_1$ and $\psi^\prime_2$ \emph{always}
have real values on the central elements $\pm z^\prime=\pm x^\prime_1x^\prime_3\cdots x^\prime_{2k-1}$
as they should--since according to Lemma \ref{centerES} the center of $\bE_{2k-1}^{1}$ is
always isomorphic to $\Z_2\times\Z_2$.

\section{Applications} In this section we describe the abstract structure of the groups $G_n$ and $H_n$
and decompose the representation $\pi_n:\B_n\rightarrow (\C^2)^{\otimes n}$ into its irreducible constituents.
   We then extend these ideas to the re-normalized representation of $\B_n$ that factors over
   the Temperley-Lieb algebra.
\subsection{$H_n$ and $G_n$ as abstract groups}
\begin{theorem}\label{Thmodd}
$H_n\cong \bE_{n-1}^{-1}$.
\end{theorem}
\begin{proof}
To verify that the map $\phi:\bE_{n-1}^{-1}\rightarrow H_n$ defined by $x_i\rightarrow g_i$
extends to a (surjective) group homomorphism one just checks that the $g_i$
satisfy the defining relations of $\bE_{n-1}^{-1}$.
Since $\Ker(\phi)$ is normal it must be trivial or intersect $Z(\bE_m^{-1})$ by Lemma \ref{centerES}(d).
We check that $\phi(-1)=\phi(x_1^2)=g_1^2=-I_{2^n}$ so $-1\not\in\Ker(\phi)$
and we have proved the theorem for $n-1$ even. If $n-1=2k-1$ is odd,
we must also check that $\pm z\not\in \Ker(\phi)$ where $z$ is the nontrivial
central element defined in Lemma \ref{centerES}.
For this we must use Lemma \ref{comps}(k) which shows that there is
a change of basis which diagonalizes the odd-indexed $g_{2i+1}$ while fixing the even indexed $g_{2i}$.
We compute the image of $z$ in this basis:
\begin{eqnarray*}
(P_n)^{-1}\phi(\pm z)P_n=(P_n)^{-1}(\pm g_1g_3\cdots g_{2k-1})P_n=\pm(\sqrt{-1})^k(\sigma_z^{\otimes 2k})
\end{eqnarray*}
which is a diagonal matrix of trace $0$, so not the identity.
\end{proof}
Combining with Theorem \ref{quotient} we have:
\begin{theorem}
The image of $\B_n$ under the representation $\pi_n$ is
an extension of $\bE_{n-1}^{-1}$ by $S_n$.\end{theorem}

\subsection{Decomposition of $\pi_n$}

By Theorem \ref{Thmodd} we have $\bE_{n-1}^{-1}\cong H_n$ as an abstract group
so the (defining) representation $(\pi_n,(\C^2)^{\otimes n})$ of $H_n$ induces a representation $\phi_n:=\pi_n\circ\phi$ of $\bE_{n-1}^{-1}$.

\subsubsection{$n$ odd}
Assume that $n=2k+1$ is odd.  Then we may decompose $(\C^2)^{\otimes 2k+1}\cong\bigoplus_i m_iV_i$
as representations of $\bE_{2k}^{-1}$ for some multiplicities $m_i$.
Let $\chi$ be the character of $\phi_{2k+1}$.  Since
$$\phi_{2k+1}(-1)=(I_2\otimes\cdots\otimes g_i^2\cdots\otimes I_2)=-I_{2^n}$$
we see that $\chi(-1)=-2^{2k+1}$ and $\chi(1)=2^{2k+1}$.
By Proposition \ref{chars} we can compute the multiplicities $m_i$ of the irreducible components $V_i$:
$$m_i=\frac{1}{2^{2k+1}}\sum_{x\in\bE_{2k}^{-1}}\chi_i(x)\overline{\chi(x)}$$
The character $\chi_1$ of the $2^k$-dimensional representation $V_1$ vanishes on the non-central
elements of $\bE_{2k}^{-1}$ so we compute the multiplicity

$$m_1=\frac{2^k\cdot 2^{2k+1}+2^k\cdot 2^{2k+1}}2^{2k+1}=2^{k+1}$$

so $V_1$ appears $2^{k+1}$ times.
But $\dim V_1=2^k$ so $\dim (2^{k+1}V_1)=2^{2k+1}=\dim(\C^2)^{\otimes 2k+1}$,
so in fact $\pi_{2k+1}$ decomposes diagonally as $2^{k+1}$ copies of the
unique $2^k$-dimensional representation $(\rho_1,V_1)$ of $\bE_{2k}^{-1}$.

\subsubsection{$n$ even}
Suppose $n=2k$ is even.  We have already established (see \ref{resV} in Section \ref{esk}) that
the restriction of the irreducible $2^k$-dimensional
representation $V_1$ of $\bE_{2k}^{-1}$ to $\bE_{2k-1}^{-1}$ decomposes as
the direct sum $W_1\oplus W_2$ of the two inequivalent irreducible $2^{k-1}$
dimensional representations $W_1$ and $W_1$.  So the $2^{2k}$-dimensional representation $\phi_{2k}$
decomposes diagonally as the direct sum of $2^{k}$ copies of each of $(\la_1,W_1)$ and $(\la_2,W_2)$.
One could also use the characters $\psi_i$ to determine these multiplicities.

\begin{rmk}
As $\pi_n(\PP_n)=\phi_n(\bE_{n-1}^{-1})$, all of the arguments above hold \emph{mutatis mutandis}
for decomposing $\pi_n$ restricted to $\PP_n$.
\end{rmk}

\subsubsection{Extension to $\B_n$}
With the explicit formulas for the representations $\rho_1$, $\la_1$ and $\la_2$ in hand,
we easily compute the extensions $\Hat{\rho}_1$, $\Hat{\la}_1$ and $\Hat{\la}_2$ to $\B_n$ using Lemma \ref{comps}(a).  Using the matrices $d$, $M$ and $D$ from Section \ref{defs} we give the explicit
matrices for the $2^k$-dimensional irreducible representation $\Hat{\rho}_1$ with $n=2k+1$ noting that the $\Hat{\la}_1\oplus\Hat{\la}_2$ is just the restriction of $\Hat{\rho}_1$.
\begin{eqnarray*}
&&\Hat{\rho}_1(\sigma_1)=d\otimes I_2^{\otimes k-1}\\
&&\quad \vdots\\
&&\Hat{\rho}_1(\sigma_{2i})=I_2^{\otimes i-1}\otimes M\otimes I_2^{\otimes k-i}\\
&&\Hat{\rho}_1(\sigma_{2i+1})=I_2^{\otimes i-1}\otimes D\otimes I_2^{\otimes k-i-1}\\
&&\quad \vdots\\
&&\Hat{\rho}_1(\sigma_{2k})=I_2^{\otimes k-1}\otimes M
\end{eqnarray*}
 The decomposition of $\pi_n$ remains the same, so
summarizing we have:
\begin{theorem}
The representation $\pi_n$ of $\B_n$ decomposes as
$$(\C^2)^{\otimes n}\cong\begin{cases} (\C^2)^{\otimes (n+1)/2}\otimes V_1 & \textrm{$n$ odd}\\
(\C^2)^{\otimes n/2}\otimes (W_1\oplus W_2) & \textrm{$n$ even}\end{cases}$$
\end{theorem}

\section{Jones representation and Jones polynomial}

The Jones representation of the braid groups $\B_n$ are defined using the
Temperley-Lieb algebras $TL_n(q)$.  Jones representation $\rho_r$
in the following means the unitary representation
 of the braid groups at $q=e^{2\pi i/r}$ factoring through the semisimple
 Temperley-Lieb algebras, which are quotients of the
 Hecke algebras in \cite{Jones86}.  The specific
 formulas that we use are the ones in \cite{FLW}.

\begin{definition}
Let $q=\sqrt{-1}$.  The Temperley-Lieb algebra $TL_n(q)$
is defined as the (semisimple) quotient of the braid group algebra $\C[\B_n]$ by
(the ideal generated by) the relations:
\begin{enumerate}
\item[TL1:] $(\sigma_i+1)(\sigma_i-q)=0$
\item[TL2:] $\sigma_i\sigma_{i+1}\sigma_i+\sigma_i\sigma_{i+1}+\sigma_{i+1}\sigma_i+\sigma_i+\sigma_{i+1}+1=0$
\item[TL3:] $(\sigma_i-\sigma_{i+1})^2=\sqrt{-1}$ (\emph{i.e.} Jones-Wenzl projector $p_3=0$)
\end{enumerate}
\end{definition}

Observing that the Yang-Baxter operator $R$ satisfies $(R-\zeta I_4)(R-\overline{\zeta}I_4)=0$
we can define a new matrix $R^\prime=-\overline{\zeta}R$ that
satisfies $(R^\prime-\sqrt{-1}I_4)(R^\prime+I_4)=0$.
Since the equation (YBE) is homogeneous, (YBE) is satisfied by $R^\prime$ also.
It is a (mildly tedious) computation to verify that the matrices
$A_1=(R^\prime\otimes I_2)$ and $A_2=(I_2\otimes R^\prime)$
satisfy $A_1A_2A_1+A_1A_2+A_2A_1+A_1+A_2+I_4=0$, and $(A_1-A_2)^2=\sqrt{-1}I_4$.
Thus the representation $\pi_n^\prime$ of $\B_n$ afforded us by $R^\prime$
(or $\C\B_n$ if we prefer) factors over the Temperley-Lieb algebra $TL_n(\sqrt{-1})$.
We can easily extend what we have learned about the representation $\pi_n$ of $\B_n$ to
this slight variation by observing the effect of renormalizing $R$.
We record the result in the following (compare to \cite{Jones86}):

\begin{cor}
Denote by $H_n^\prime=\pi_n^\prime(\PP_n)$ and $G_n^\prime=\pi_n^\prime(\B_n)$.
Then we have $H_n^\prime\cong\bE_{n-1}^1$, and $G_n^\prime/H_n^\prime\cong S_n$.
\end{cor}
\begin{proof} This follows easily from the observation that renormalizing $R$
by $-\overline{\zeta}$ has the effect of multiplying
the generators $g_i$ of $H_n$ by $-\sqrt{-1}$.
Doing the same to the generators of the group $\bE_{n-1}^{-1}$
just gives us a presentation of the group $\bE_{n-1}^1$,
and the same arguments as in the original representation $\pi_n$ go through verbatim.
\end{proof}

To relate $\pi_n^\prime$ to the Jones representation $\rho_4$ of $\B_n$,
we recall some facts about the
Jones representation.  The Temperley-Lieb algebras at a 4-th root
of unity are
complex Clifford algebras and are isomorphic to
the matrix algebra of $2^{n-1}\times 2^{n-1}$ matrices if $n$ is odd, and the
direct sum of two matrix algebras of $2^{\frac{n}{2}-1}\times 2^{\frac{n}{2}-1}$
matrices if $n$ is even \cite{Jones}.
(Note here $n$ is
the number of strands in the geometric realization of $\B_n$, and differs by 1 from Jones'
notation in \cite{Jones86}.)  So the Jones representation $\rho_4$ consists of
a single irreducible sector if $n$ is odd, and the
direct sum of two irreducible sectors if $n$ is even.
 Comparing with the comments in Subsection \ref{nuis1reps} we can also
determine the decomposition of the representation $\pi_n^\prime$ as before.
 It follows that the
restriction of the Jones representation $\rho_4$ to $\PP_n$ for $n$ even is the
odd representation $V_1$ of the extra-special 2-group $\bE_{n-1}^1$, and
for $n$ odd, $W_1\oplus W_2$ as in Theorem 4.4.  The images $\rho_4(\B_n)$ fit
into the following exact sequence:
$$ 1\rightarrow \bE_{n-1}^1 \rightarrow \rho_4(\B_n)\rightarrow
S_n\rightarrow 1.$$
Projectively, we have
$$ 1\rightarrow \Z_2^{n-1}\rightarrow \rho_4(\B_n)\rightarrow
S_n\rightarrow 1.$$
The symmetric group $S_n$ acts on the
coordinates of  $\Z_2^{n}$, hence $\Z_2^{n-1}$ when $n$ is even.
This action splits the exact sequence.
But when $n$ is odd, this sequence does not split as is shown in \cite{Jones86}.

The Jones polynomial of a link at $\sqrt{-1}$ is given by the following formula \cite{FLW}:
$$J_4(\hat{\sigma})=(-1)^{n-1+\frac{e(\sigma)}{4}}\cdot
(\sqrt{2})^{-\frac{1+(-1)^n}{2}}\cdot \textrm{Trace}(\rho_4(\sigma)),$$
where $e(\sigma)$ is the sum of all exponents of the standard
braid generators appearing in $\sigma$, and $\hat{\sigma}$ is the closure of $\sigma$.
We can also define link invariants using the flipped R-matrix $R$.  The
conditions for enhancement $(\mu_i, \alpha, \beta)$ is given in \cite{Turave} Theorem 2.3.1.
Working through
the conditions, we found two link invariants:
$T_R(\hat{\sigma},\alpha)=\alpha^{n-e(\sigma)}\cdot (\sqrt{2})^{-n}\cdot
\textrm{Trace}(\pi_n(\sigma))$,
where $\alpha=\pm 1$.  Comparing with the Jones polynomial we get the
relation:
$$T_R(\hat{\sigma},\alpha)=(-1)^{n-1+e(\sigma)}\cdot {\alpha}^{n-e(\sigma)}
\cdot \sqrt{2}\cdot J_4(\hat{\sigma}).$$
As we know that Jones polynomial $J_4(\hat\sigma)$ is
$(\sqrt{2})^{c(\hat{\sigma})-1}\cdot (-1)^{\textrm{Arf}(\hat{\sigma})}$
if $\textrm{Arf}(\hat{\sigma})$ is defined and $0$ otherwise, where
$c(\hat{\sigma})$ is the number of components of the link $\hat{\sigma}$ \cite{Jones}.
It follows that
$T_R(\hat{\sigma},\alpha)$ computes essentially the Arf invariant
of a link.

\end{document}